# Random Consolidations and Fragmentations Cycles Lead to Benford's Law

**ABSTRACT**

Benford's Law predicts that the first significant digit on the leftmost side of numbers in real-life data is proportioned between all possible 1 to 9 digits approximately as in LOG(1 + 1/digit), so that low digits occur much more frequently than high digits in the first place. For example, digit 1 occurs approximately 30.1% in the first place in random numbers, while digit 9 occurs only approximately 4.6%. In this article it is shown that a process where a large enough set of identical quantities constantly alternates between minuscule random consolidations (summing two randomly chosen values into a singular value) and tiny random fragmentations (division of one randomly chosen value into two new values) converges digit-wise to the Benford proportions after sufficiently many such cycles. The statistical tendency of the system after numerous cycles is to have approximately 2/3 multiplicative expressions which are conducive to Benford behavior as they tend to the Lognormal Distribution, and 1/3 additive expressions which are detrimental to Benford behavior as they tend to the Normal Distribution, hence the process represents in essence a tug of war between addition and multiplication. Since the process encounters the so-called Achilles' heel of the Central Limit Theorem, namely additions of skewed distributions with high order of magnitude, additions are not very effective, and the war is decisively won by multiplication, leading to Benford behavior. Randomness in selecting the particular quantity to be fragmented, as well as randomness in selecting the two particular quantities to be consolidated, is essential for convergence. Not surprisingly then, fragmentation itself could be performed either randomly say via a realization from the continuous Uniform on (0, 1), or deterministically via any fixed split ratio such as say 25% - 75%, and Benford's Law emerges in either case.

Alex Ely Kossovsky
akossovsky@gmail.com

# PART 1: BENFORD'S LAW

# PART 2: RANDOM CONSOLIDATIONS AND FRAGMENTATIONS PROCESSES

# PART 1:
# BENFORD'S LAW

# [I]   The First Digit on the Left Side of Numbers

It has been discovered that the first digit on the left-most side of numbers in real-life data sets is most commonly of low value such as {1, 2, 3} and rarely of high value such as {7, 8, 9}. As an example serving as a brief and informal empirical test, a small sample of 40 values relating to geological data on time between earthquakes is randomly chosen from the data set on all global earthquake occurrences in 2012 – in units of seconds. Figure A depicts this small sample of 40 numbers. Figure B emphasizes in bold and black color the 1st digits of these 40 numbers.

| 285.29  | 185.35  | 2579.80 | 27.11   |
|---------|---------|---------|---------|
| 5330.22 | 1504.49 | 1764.41 | 574.46  |
| 1722.16 | 815.06  | 3686.84 | 1501.61 |
| 494.17  | 362.48  | 1388.13 | 1817.27 |
| 3516.80 | 5049.66 | 2414.06 | 387.78  |
| 4385.23 | 2443.98 | 2204.12 | 1224.42 |
| 1965.46 | 3.61    | 1347.30 | 271.23  |
| 3247.99 | 753.80  | 1781.45 | 593.59  |
| 1482.64 | 1165.04 | 4647.39 | 1219.19 |
| 251.12  | 7345.52 | 1368.79 | 4112.13 |

**Figure A**: Sample of 40 Time Intervals between Earthquakes

| **2**85.29  | **1**85.35  | **2**579.80 | **2**7.11   |
|-------------|-------------|-------------|-------------|
| **5**330.22 | **1**504.49 | **1**764.41 | **5**74.46  |
| **1**722.16 | **8**15.06  | **3**686.84 | **1**501.61 |
| **4**94.17  | **3**62.48  | **1**388.13 | **1**817.27 |
| **3**516.80 | **5**049.66 | **2**414.06 | **3**87.78  |
| **4**385.23 | **2**443.98 | **2**204.12 | **1**224.42 |
| **1**965.46 | **3**.61    | **1**347.30 | **2**71.23  |
| **3**247.99 | **7**53.80  | **1**781.45 | **5**93.59  |
| **1**482.64 | **1**165.04 | **4**647.39 | **1**219.19 |
| **2**51.12  | **7**345.52 | **1**368.79 | **4**112.13 |

**Figure B**: The First Digits of the Earthquake Sample

Clearly, for this very small sample, low digits occur by far more frequently on the first position than do high digits. A summary of the digital configuration of the sample is given as follows:

Digit Index: { 1, 2, 3, 4, 5, 6, 7, 8, 9 }
Digits Count totaling 40 values: { 15, 8, 6, 4, 4, 0, 2, 1, 0 }
Proportions of Digits with '%' sign omitted: {38, 20, 15, 10, 10, 0, 5, 3, 0 }

Assuming (correctly) that these 40 values were collected in a truly random fashion from the large data set of all 19,452 earthquakes occurrences in 2012; without any bias or attempt to influence first digits occurrences; and that this pattern is generally found in many other data sets, one then may conclude with the phrase "not all digits are created equal", or rather "not all <u>first</u> digits are created equal", even though this seems to be contrary to intuition and against all common sense.

The focus here is actually on the <u>first meaningful digit</u> – counting from the left side of numbers, excluding any possible encounters of zero digits which only signify ignored exponents in the relevant set of powers of ten of our number system. Therefore, the complete definition of the **First Leading Digit** is the <u>first non-zero digit</u> of any given number on its left-most side. This digit is the first significant one in the number as focus moves from the left-most position towards the right, encountering the first non-zero digit signifying some quantity; hence it is also called the **First Significant Digit**. For 2365 the first leading digit is 2. For 0.00913 the first leading digit is 9 and the zeros are discarded; hence even though strictly-speaking the first digit on the left-most side of 0.00913 is 0, yet, the first significant digit is 9. For the lone integer 8 the leading digit is simply 8. For negative numbers the negative sign is discarded, hence for -715.9 the leading digit is 7. Here are some more illustrative examples:

**6**,719,525 → digit 6
0.0000**7**61 → digit 7
-0.**2**81264 → digit 2
**8**75 → digit 8
**3** → digit 3
-**5** → digit 5

For a data set where all the values are greater than or equal to 1, such as in the sample of the earthquaqe data, the first digit on the left-most side of numbers is also the First Leading Digit and the First Significant Digit, and necesarily one of the nine digits {1, 2, 3, 4, 5, 6, 7, 8, 9}; while digit 0 never occurs first on the left-most side.

# [II]   Benford's Law and the Predominance of Low Digits

Benford's Law states that:

Probability[First Leading Digit is d] = $LOG_{10}(1 + 1/d)$

$LOG_{10}(1 + 1/1) = LOG(2.00) = 0.301$
$LOG_{10}(1 + 1/2) = LOG(1.50) = 0.176$
$LOG_{10}(1 + 1/3) = LOG(1.33) = 0.125$
$LOG_{10}(1 + 1/4) = LOG(1.25) = 0.097$
$LOG_{10}(1 + 1/5) = LOG(1.20) = 0.079$
$LOG_{10}(1 + 1/6) = LOG(1.17) = 0.067$
$LOG_{10}(1 + 1/7) = LOG(1.14) = 0.058$
$LOG_{10}(1 + 1/8) = LOG(1.13) = 0.051$
$LOG_{10}(1 + 1/9) = LOG(1.11) = 0.046$
                                    ---------
                                      1.000

Figure C depicts the distribution. Figure D visually depicts Benford's Law as a bar chart. This set of nine proportions of Benford's Law is sometimes referred to in the literature as **'The Logarithmic Distribution'**. Remarkably, Benford's Law is confirmed in almost all real-life data sets with high order of magnitude, such as in data relating to physics, chemistry, astronomy, economics, finance, accounting, geology, biology, engineering, governmental census data, and many others.

| Digit | Probability |
|---|---|
| 1 | 30.1% |
| 2 | 17.6% |
| 3 | 12.5% |
| 4 | 9.7% |
| 5 | 7.9% |
| 6 | 6.7% |
| 7 | 5.8% |
| 8 | 5.1% |
| 9 | 4.6% |

**Figure C**:  Benford's Law for First Digits

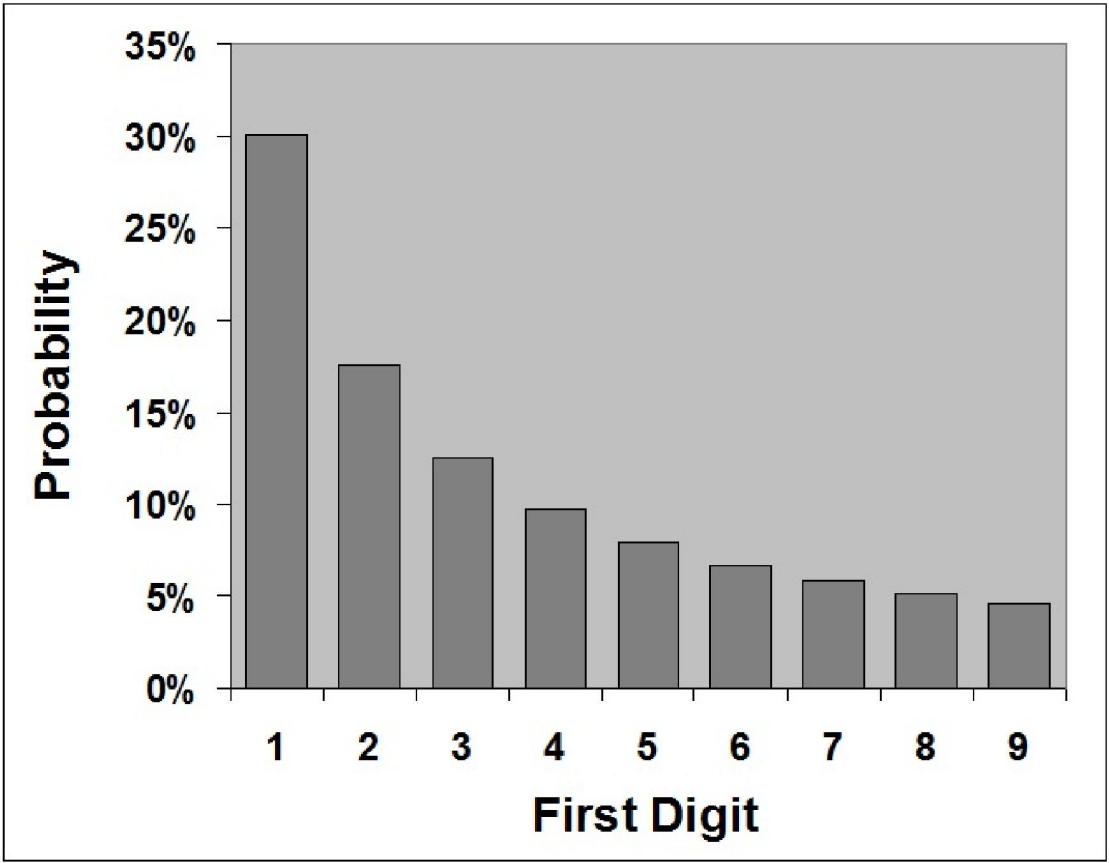

**Figure D**: Benford's Law – Probability of First Leading Digit Occurrences as a Bar Chart

# [III] Sum of Squares Deviation Measure (SSD)

It is necessary to establish a standard measure of 'distance' from the Benford digital configuration for any given data set. Such a numerical measure could perhaps tell us about the conformance or divergence from the Benford digital configuration of the data set under consideration. This is accomplished with what is called **Sum Squares Deviations (SSD)** defined as the sum of the squares of the 'errors' between the Benford expectations and the actual/observed values (in percent format – as opposed to fractional/proportional format):

$$SSD = \sum_{1}^{9} \left( \text{Observed \% of d} - 100 * \text{LOG}\left(1 + \frac{1}{d}\right) \right)^2$$

For example, for observed 1st digits proportions of {31.1, 18.2, 13.3, 9.4, 7.2, 6.3, 5.9, 4.5, 4.1} with '%' sign omitted, SSD measure of distance from the logarithmic is calculated as:

**SSD** = $(31.1 – \textbf{30.1})^2 + (18.2 – \textbf{17.6})^2 + (13.3 – \textbf{12.5})^2 + (9.4 – \textbf{9.7})^2 +$
$+ (7.2 – \textbf{7.9})^2 + (6.3 – \textbf{6.7})^2 + (5.9 – \textbf{5.8})^2 + (4.5 – \textbf{5.1})^2 + (4.1 – \textbf{4.6})^2 = \textbf{3.4}$

SSD generally should be below 25; a data set with SSD over 100 is considered to deviate too much from Benford; and a reading below 2 is considered to be ideally Benford.

# [IV]  Physical Order of Magnitude of Data (POM)

Rules regarding expectations of compliance with Benford's Law rely heavily on measures of order of magnitude (i.e. variability) of data.

Physical order of magnitude of a given data set is a measure that expresses the extent of its variability. It is defined as the ratio of the maximum value to the minimum value. The data set is assumed to contain only positive numbers greater than zero.

Physical Order of Magnitude (POM) = Maximum/Minimum

The classic definition of order of magnitude involves also the application of the logarithm to the ratio maximum/minimum, transforming it into a smaller and more manageable number.

Order of Magnitude (OOM) = $LOG_{10}$(Maximum/Minimum)
Order of Magnitude (OOM) = $LOG_{10}$(Maximum) - $LOG_{10}$(Minimum)

Since such logarithmic transformation has a monotonic one-to-one relationship with max/min, it does not provide for any new insight or information, but could rather be looked upon in a sense as simply the use of an alternative scale, still measuring the same thing. For this reason the complexity of logarithm can be avoided altogether by referring only to the simple POM measure.

The more profound reason for using **POM** instead of **OOM** is its feature as a universal measure of variability, totally independent of societal number system in use, as well as being independent on the arbitrary choice of base 10, derived from the chanced or random occurrence of us having 10 fingers. This is the motivation behind the use of the term 'physical', expressing real and physical measure of variability, divorced from any numerical inventions, and especially so when data relates to the natural world such as in scientific figures and physical information.

OOM is perhaps more appropriate for a single isolated number, where it is re-defined as simply $LOG_{10}$(Number) without any reference to maximum, minimum, or any ratio. If we can assume that that number is an integral whole number without any [trailing] fractional part, then an alternative meaning of this OOM definition is simply expressing how many digits approximately are necessary to write the number. Surely in general, the bigger the [integral] number the more digits it takes to write it! For example, $LOG_{10}$(8,200,135) = 6.9, which is about 7, and that's exactly how many digits the number involves. As another example, $LOG_{10}$(10,000,000) = 7.0 which is exactly one digit less than the number of digits involved in writing the number, namely 8 digits.

In the extreme case where all the numbers in the data set are identical, having the value R say, variability is then nonexistent, and POM = maximum/minimum =  R/R = 1.

NOTE: LOG(X) or log(x) notation in this article would always refer to our decimal base 10 number system, hence the more detailed notation of $LOG_{10}(X)$ is often (but not always) avoided.

# [V]  A Robust Measure of Physical Order of Magnitude (CPOM)

It is perhaps unfortunate that the literature in statistics does not seem to contain any robust definition of order of magnitude. Such a measure should prove steady and consistent for all types of data sets, strongly resisting outliers, preventing them from overly influencing the numerical measure of data variability.

In order to accomplish exactly that, and also to preserve the advantage of avoiding dependencies on arbitrary societal number systems and particular bases, the basic (independent) structure of POM shall be used, but with the added modification of simply eliminating any possible outliers on the left for small values and on the right for big values. This is accomplished by narrowing the focus exclusively onto the core 80% part of the data. This brutal purge eliminates any malicious and misleading outliers as well as any innocent and proper data points which happened to stray just a little bit away from the core part of the data. The measure shall be called Core Physical Order of Magnitude and it is defined as follows:

Core Physical Order of Magnitude (CPOM) = $Q_{90\%} / Q_{10\%}$

The definition simply reformulates POM by substituting the 10th percentile (in symbols $\mathbf{Q_{10\%}}$) for the minimum, and by substituting the 90th percentile (in symbols $\mathbf{Q_{90\%}}$) for the maximum.

The 10th percentile is the value below which about 10% of the data points may be found. The 50th percentile is the median, below which about half of all the data points may be found. The 90th percentile is the value below which about 90% of the data points may be found.

In general, the rejection of outliers appearing in a given data set may be justified, or it may actually be misguided. For example, if over 50,000 students at a large university are surveyed with regards to height, and the top value is say 7.25 meter, then this outlier is certainly some kind of an error in recording and should be excluded from further analysis. If the top value is say 2.37 meter, then this 2.37 outlier is actually an integral part of the data set, and especially so if the well-known tall student is ordered to appear at the administration office, rudely interrupting his exciting basketball game at the court, and another measurement is taken, confirming his 2.37 meter height as well as his existence. This is not simply a matter of mere semantics, and there is a compelling argument not to classify this 2.37 value as an outlier, although in reality it depends on the context.

# [VI]  Two Essential Requirements for Benford Behavior

One of the two essential prerequisites or conditions for data configuration with regards to compliance with Benford's Law is that the value of the order of magnitude of the data set should be approximately over 3; in other words, that $LOG_{10}$(Maximum/Minimum) > 3, and that therefore (Maximum/Minimum) > $10^3$. This in turn implies that the threshold POM value (separating compliance from non-compliance) is about 1000, namely that POM > 1000 constitutes the condition for compliance.

The above prerequisite for compliance totally ignores the thorny issue of outliers and edges, and in that sense it is too simplistic and even completely erroneous for some data sets. Hence, using the CPOM qualification is essential in judging whether or not a given data set complies or does not comply with Benford's Law. The proper qualification for compliance with the law in the approximate - obtained via extensive empirical studies - is then as follows:

Core Physical Order of Magnitude  = $Q_{90\%}$ / $Q_{10\%}$ > 100

Actually, even lower CPOM values such as 50 and 30 are Benford, but falling below 30 does not bode well for getting anywhere near the logarithmic distribution.

Skewness of data where the histogram comes with a prominent tail falling to the right is the second essential criterion necessary for Benford behavior. Indeed, most real-life physical data sets are generally skewed in the aggregate, so that overall their histograms have tails falling on the right, and consequently the quantitative configuration is such that the small is numerous and the big is rare, while low first digits decisively outnumber high first digits.

The <u>asymmetrical</u>, Exponential, Lognormal, k/x [and many other distributions] are typical examples of such quantitatively skewed configuration, and therefore they are approximately, nearly, or exactly Benford - respectively. The <u>symmetrical</u> Uniform, Normal, Triangular, Circular-like, and other such distributions are inherently non-Benford, or rather anti-Benford, as they lack skewness and do not exhibit any bias or preference towards the small and the low.

Symmetrical distributions are always non-Benford, no matter what values are assigned to their parameters. By definition they lack that asymmetrical tail falling to the right, and such lack of skewness precludes Benford behavior regardless of the value of their order of magnitude.  Order of magnitude simply does not play any role whatsoever in Benford behavior for symmetrical distributions. For example, first digits of the Normal($10^{35}$, $10^8$) or the Uniform(1, $10^{27}$) are not Benford at all, and this is so in spite of their extremely large orders of magnitude. In summary: Benford behavior in extreme generality can be found with the confluence of sufficiently large order of magnitude together with skewness of data - having a histogram falling to the right. The combination of skewness and large order of magnitude is not a guarantee of Benford behavior, but it is a strong indication of likely Benford behavior under the right conditions. Moderate [overall] quantitative skewness with a tail falling too gently to the right implies that digits are not as skewed as in the Benford configuration. Extreme [overall] quantitative skewness with a tail falling sharply to the right implies that digits are severely skewed, even more so than they are in the Benford configuration.

Bowley Skewness for example, defined as [(Q3 – Q2) – (Q2 – Q1)] **/** [Q3 – Q1] is an intuitive measure of skewness but its numerical value fluctuates greatly across data sets. Calculated Bowley Skewness values for numerous logarithmic data sets and distributions do not yield any consistent result, except that all values come out above 0.3 and below 1.0, and which is consistent with the fact that all logarithmic data sets are positively skewed in the aggregate. In sharp contrast, almost all non-logarithmic data sets come out with decisively lower Bowley Skewness values below 0.25 and above 0. In contrast to Bowley's unstable value for logarithmic data sets, Benford's Law is a very consistent and almost exact measure of skewness, with very little fluctuations across logarithmic data sets.

**Related Log Conjecture** in Benford's Law states that whenever the density curve of the logarithm of the data starts from the bottom on the log-axis itself, rises up to some plateau or maximum point, and then falls back all the way down to the log-axis itself, mimicking in a sense an upside-down U-like curve, and where total span on the log-axis from the minimum log to the maximum log is at least 3 units, then these two essential prerequisites above are guaranteed to be satisfied, and in addition guaranteeing that the data itself is nearly or almost perfectly Benford. See Kossovsky (2014) chapter 63 for a full discussion and concrete examples about the conjecture.

## [VII]   Data Skewness is More Prevalent than Benford's Law

All data sets obeying Benford's Law (i.e. logarithmic data) are structured in such a way that in the aggregate there are more small quantities than big quantities. In other words, that in the aggregate the histogram is falling to the right, except perhaps in the beginning on the very left for low values where it temporarily rises for a very small portion of overall data, as well as in few and minors reversals along the way where it rises briefly. This quantitative configuration is called 'positive skewness' in mathematical statistics.

This highly prevalent phenomenon of skewness has by far much wider scope and it is much more common in the physical world and in the realm of abstract mathematics than the more particular Benford quantitative configuration. This statement does not imply that Benford's Law is not prevalent in scientific, physical, and numerous other data types, on the contrary, it is highly prevalent. The statement only implies that in almost all the counter examples and exceptions to Benford's Law, the phenomenon of skewness is still found, albeit with different quantitative configurations than that of the Benford one (and typically with milder skewness, but at times even skewer). This renders Benford's Law a subset of the more universal skewness phenomenon. The assertion is derived from concrete experience with numerous real-life numerical examples and from general research in Benford's Law. While this discussion may sound vague, in fact it is rather a very essential overview in the entire quantitative phenomenon of Benford's Law and of real-life data analysis in general. For those statisticians and data analysts who have worked on data sets and the Benford phenomenon for many years, including doing theoretical research, this generic statement seems natural, fundamental, and quite necessary.

# PART 2:

# RANDOM CONSOLIDATIONS AND FRAGMENTATIONS PROCESSES

# [1]  Random Consolidations and Fragmentations Cycles Lead to Benford

In this article it is shown that a process where a large set of identical quantities constantly alternates between minuscule random consolidations (summing two randomly chosen values into a singular value) and tiny random fragmentations (division of one randomly chosen value into two new values) converges eventually to the Benford configuration after sufficiently many cycles. Randomness in selecting the particular quantity to be fragmented, as well as randomness in selecting the two particular quantities to be consolidated, is essential for convergence. Not surprisingly then, fragmentation itself could be performed either randomly say via a realization from the continuous Uniform on (0, 1), or deterministically via any fixed ratio of breakup such as say 25% - 75% or 10% - 90%, and the Benford configuration emerges in either case. This is so since randomness already exists in the system in the determination of which quantities are to be fragmented and consolidated.

The random model discussed in this article is a cyclical process which alternates constantly between (1) the selection of one quantity chosen at random and its breakdown [fragmentation] into two new smaller quantities, (2) followed by the selection of two randomly chosen quantities and their summation [consolidation] into a singular added quantity. Monte Carlo computer simulations results decisively show that such a process leads to a rapid Benford digital convergence. Exploration of the resultant (random) algebraic expressions coming out of the corresponding mathematical model, nicely explains this Benford behavior simply in terms of a random multiplication process, and this is so in spite of the addition terms involved here. In other words, the arithmetical model of the process points to a tug of war between addition and multiplication, and where multiplication ultimately and decisively triumphs over addition.

This consolidations and fragmentations model, abbreviated as **C&F**, is based on **L** balls, all with an identical initial real value **V**, such as for example the initial weight of each ball, and assuming uniformity of mass density, implying equivalency between volume and weight proportions. Each C&F model consists of **C** cycles. Within one full cycle, two opposing processes are performed, one of fragmentation, followed by one of consolidation. The first process is the fragmentation of a single ball chosen at random utilizing the discrete Uniform on {1, 2, 3, … , L}, as well as the continuous Uniform(0, 1) to decide on the proportions of the two fragments. The second process is the consolidation of two balls chosen at random and fused into a singular and larger ball, utilizing at first the discrete Uniform distribution on {1, 2, 3, … , L, L + 1}, followed by the utilization of the discrete Uniform on {1, 2, 3, … , L}, for ball selections. Obviously, after each full cycle, the number of existing balls is unchanged, and it is still L, being the same as in the beginning and as in the end of the entire process. Also, the total quantity of the entire system - namely the overall sum or overall weight of all the balls - is conserved throughout the entire process. The specific description of physical balls of uniform mass density being broken and then fused, and the focus on the weight variable, is an arbitrary one of course, and the generic model is of pure quantities and abstract numbers.

Schematically the process is described as follows:

1) Initial Set = {V, V, V, …  L times   …V, V, V}
2) Repeat C times:
    (i) Choose one ball at random and split it as in {Uniform(0, 1), 1 – Uniform(0, 1)}.
    (ii) Choose two balls randomly and merge them.
3) Final Set is Benford.

For example, for 7 balls with an initial value of 35 each, that is L = 7, V = 35, we record the initial few cycles as the process develops randomly step-by-step:

{35, 35, 35, 35, 35, 35, 35}
{35, 35, **35**, 35, 35, 35, 35}  to be split
{35, 35, 31, 4, 35, 35, 35}
{35, 35, 31, 4, **35**, 35, **35**}  to be merged
{35, 35, 31, 4, 35, 70, 35}
{35, 35, 31, **4**, 35, 70, 35}  to be split
{35, 35, 31, 2, 2, 35, 70, 35}
{35, 35, 31, 2, **2**, 35, **70**, 35}  to be merged
{35, 35, 31, 2, 72, 35, 35}
{**35**, 35, 31, 2, 72, 35, 35}  to be split
{25, 10, 35, 31, 2, 72, 35, 35}
{25, 10, 35, 31, 2, 72, **35**, **35**}  to be merged
{25, 10, 35, 31, 2, 72, 70}

Clearly for such tiny set of balls where L = 7 there could never be any convergence to the logarithmic, even if one goes much further than C = 3 cycles. In fact, at least two digits out of all possible 9 first digits would obtain the embarrassing very low 0% proportion here no matter. It should be noted that the simulation above deliberately avoided fractional values, for pedagogical purposes, keeping the values presented as simple as possible for easy demonstrations.

Monte Carlo computer simulations show that after C full cycles, the weight of the balls is very nearly Benford, given that C > 2*L approximately, although C > 3*L or C > 4*L usually give slightly better results. It is necessary to have a sufficient number of these C cycles so that all or almost all of the balls experience either fragmentation or consolidation (preferably both, and hopefully not merely once, but rather twice or three times). By cycling at least twice as many balls that exist in the system, we ensure that (almost) all the balls undergo transformation of some sort, and that the initial value V is (almost) nowhere to be found among the balls at the end of the entire process. Continuing beyond the required 2*L or 3*L cycles does not ruin the logarithmic convergence thus obtained, and Benford is steadily preserved (or rather further perfected) as more cycles are added. Surely, the other essential prerequisite for logarithmic convergence here is to have a sufficiently large number of balls in the system so that the logarithmic can be properly manifested. Hence the requirement is that L > 200, or for better convergence that L > 300. Falling below 300 or 200 balls for example yields only crude or approximate logarithmic-like results, as in all small data sets aspiring to obey the law.

Core Physical Order of Magnitude (CPOM) shall be in use here to ensure robustness.
All C&F processes start out with V as the unique and repeated value, having no variability whatsoever, hence CPOM initially is the lowest possible value of $Q_{90\%}/Q_{10\%} = V/V = 1$. But at the end of the entire C&F process, CPOM is sufficiently high and Benford behavior is found.

## [2] Empirical Examinations of Consolidations and Fragmentations Models

We now run Monte Carlo computer simulations for three different Consolidations and Fragmentations schemes and explore results.

**Scheme A:**

In Monte Carlo computer simulation of one particular C&F process, 2000 balls, with an initial identical value of 1, undergo 8000 binary cycles, utilizing the random Uniform(0, 1) for fragmentations. The computerized results obtained here were as follows:

Initial Set - {1, 1, 1, …   2000 times   … , 1, 1, 1}
Final Set  - {0.0000000012, 0.0000000039, 0.0000000227,   …   15.5, 16.5, 16.9}

Only the first 3 and the last 3 elements from the ordered final set above are shown, namely the extreme values, the 3 biggest and the 3 smallest values, surely to be considered as outliers.

First Digits of Final Set  -  {29.1, 19.1, 11.9, 9.8, 8.3, 6.8, 6.4, 4.7, 4.2}
Benford's Law 1st Digits - {30.1, 17.6, 12.5, 9.7, 7.9, 6.7, 5.8, 5.1, 4.6}

Figure 1 depicts the digital comparison between this C&F process and Benford's Law.
The low 4.7 SSD value calculated here indicates that results are very close to the Benford configuration.

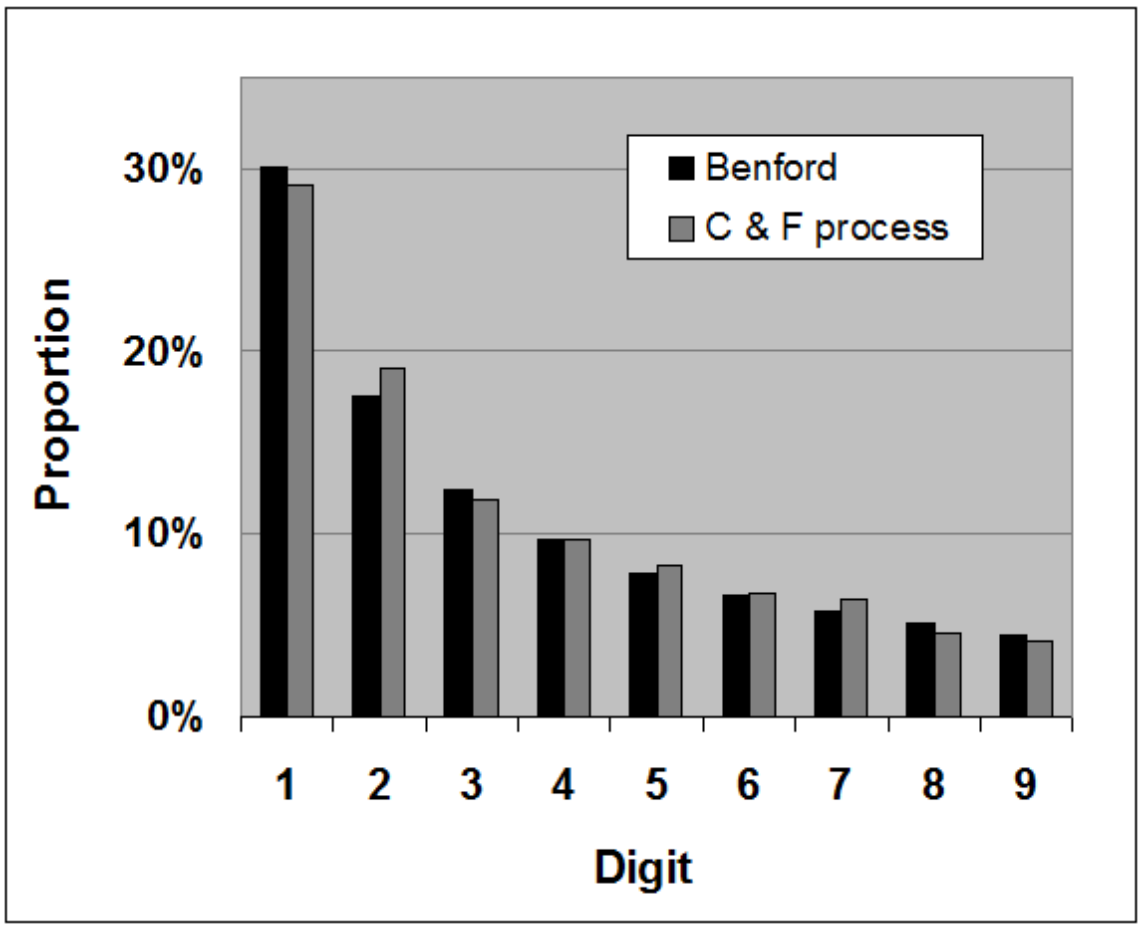

**Figure 1**: C&F Process with a Random Fragmentation Ratio Uniform(0, 1) – Scheme A

Quantitatively, the process takes as input a set of identical numbers with CPOM value of 1, and transforms it into a highly skewed set of numbers, where the small is much more numerous than the big, and where CPOM = $Q_{90\%}/Q_{10\%}$ = 3.05/0.0014 = 2131. Here 80% of resultant data falls within (0.0014, 3.05).

Figure 2 depicts the histogram from 0 to 10.2 of the final resultant set of numbers after the 8000th binary cycle. For better visual clarity a logarithmic vertical scale is used, although this masks the dramatic quantitative fall occurring here.

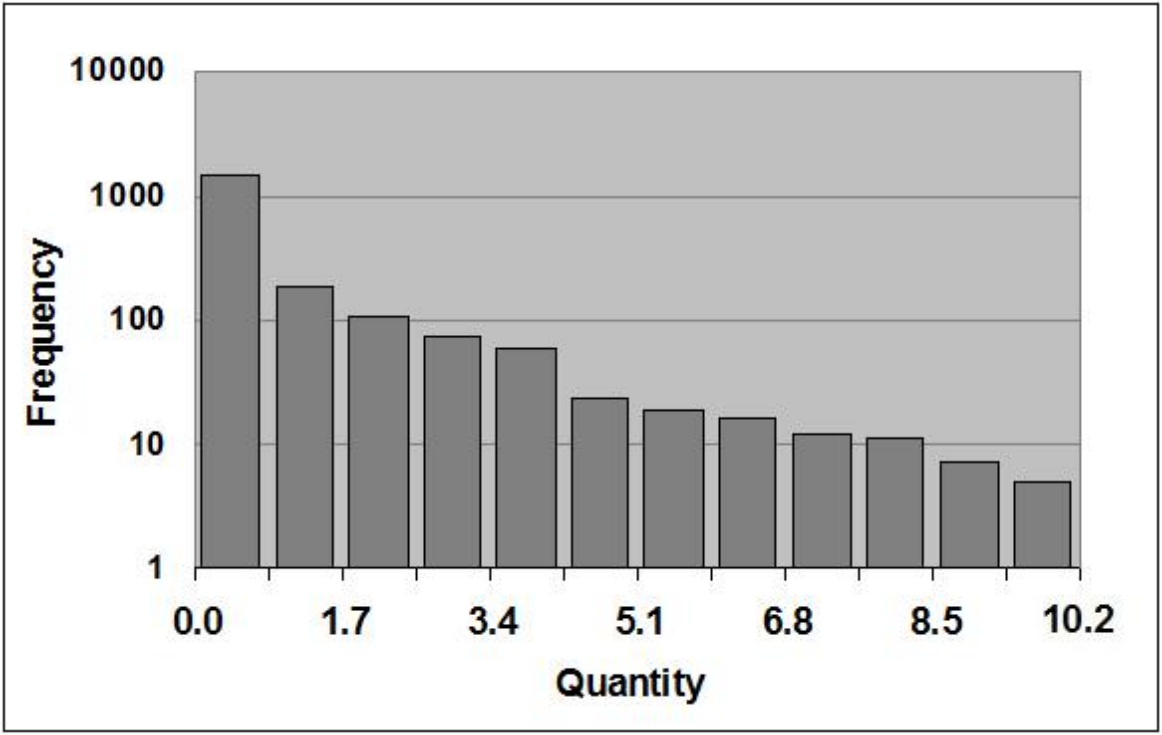

**Figure 2**: Quantitative Configuration with a Random Ratio Uniform(0, 1)   –   Scheme A

Significantly, the strong tendency towards the logarithmic in consolidations and fragmentations processes is such that randomness regarding the proportions of the fragments of any split ball is not even required. Instead of utilizing the random continuous Uniform(0, 1) to split a ball into two pieces, one may utilize a deterministic fixed ratio called p, and the process converges to Benford just as rapidly! The necessary factor driving convergence here is only the random manner by which the ball which is to be broken is chosen, and the random manner by which the two balls which are to be consolidated are chosen. Such randomness is generated here via the discrete uniform distributions [although other more complex versions and models could be devised and imagined which randomly select balls by taking their sizes into consideration]. No obvious constraint whatsoever seems to exist here for the value of p, and rapid convergence to the logarithmic is found for all the many simulation trials performed by the author with varying p values. Moreover, there seems to be no distinction whatsoever in the rate of convergence whether one utilizes a random model or a deterministic model in splitting a ball. This is so since there exists already 'plenty of randomness' in the system in how balls are selected for fragmentations and for consolidations.

**Scheme B:**

In another Monte Carlo simulation example utilizing a deterministic fixed 50% - 50% ratio, with 1500 balls all having the identical initial value of 1, being cycled 10,000 times, first digits came out nearly Benford:

Initial Set - {1, 1, 1, …     1500 times     … , 1, 1, 1}
Final Set  - {0.0000124,  0.0000153,  0.0000153,     …      13.6,  18.1,  18.6}

Only the first 3 and the last 3 elements from the ordered final set above are shown, namely the extreme values, the 3 biggest and the 3 smallest values, surely to be considered as outliers.

First Digits of Final Set   -   {31.3,  17.0,  12.9,  9.2,  6.0,  7.5,  5.7,  5.6,  4.8}
Benford's Law 1st  Digits -  {30.1,  17.6,  12.5,  9.7,  7.9,  6.7,  5.8,  5.1,  4.6}

Figure 3 depicts the digital comparison between this C&F process and Benford's Law. The low 6.7 SSD value calculated here indicates that results are very close to the Benford configuration.

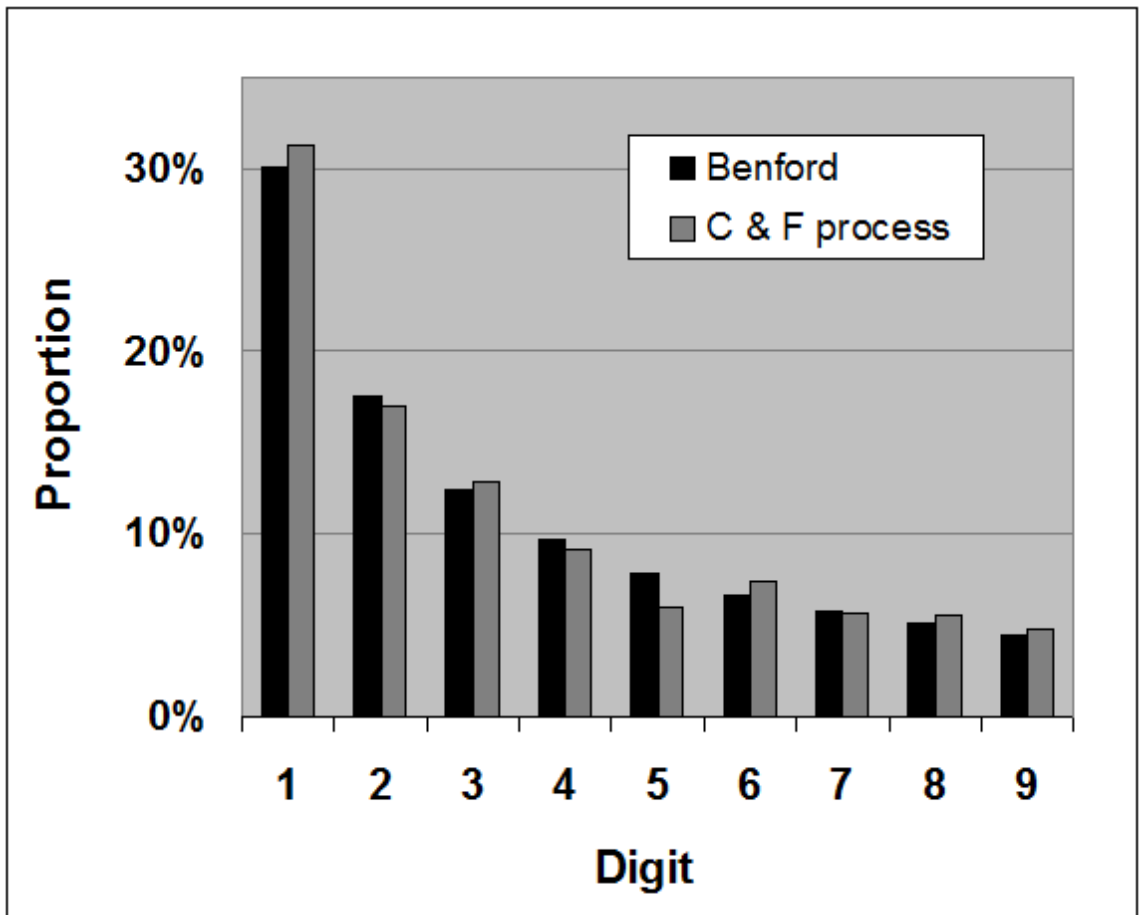

**Figure 3**: C&F Process with a Deterministic Fixed Even Ratio 50% - 50%  –   Scheme B

Quantitatively, the process takes as input a set of identical numbers with CPOM value of 1, and transforms it into a highly skewed set of numbers where the small is much more numerous than the big, and where CPOM is $Q_{90\%}/Q_{10\%} = 2.78/0.0115 = 241$. Here 80% of resultant data falls within (0.0115, 2.78). Figure 4 depicts the histogram from 0 to 10.5 of the final resultant set of numbers after the 10,000th binary cycle.

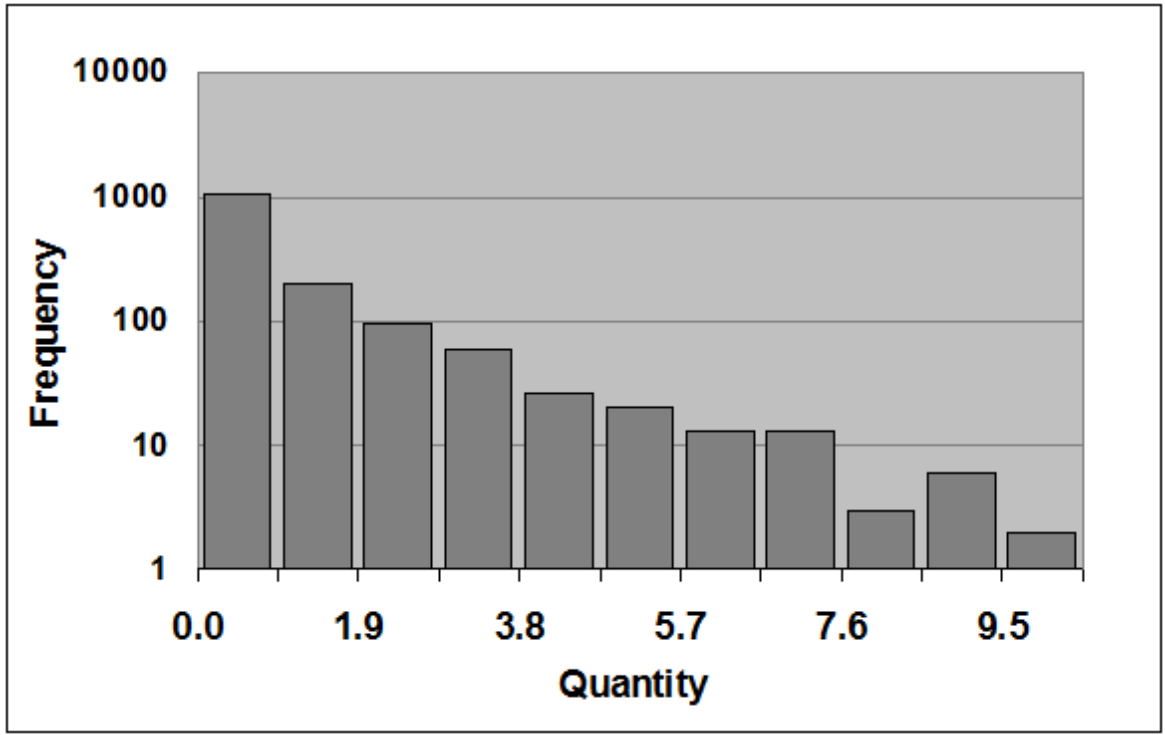

**Figure 4**: Quantitative Configuration with a Fixed Even Ratio 50% - 50% – Scheme B

**Scheme C:**

In another Monte Carlo simulation example utilizing a deterministic fixed 15% - 85% ratio, with 1000 balls all having the identical initial value 1, being cycled 3000 times, first digits came even closer to Benford:

Initial Set - {1, 1, 1, …    1000 times    … , 1, 1, 1}
Final Set  - {0.0000000278,  0.0000000278,  0.0000001574,    …    , 13.1, 13.5, 15.4}

Only the first 3 and the last 3 elements from the ordered final set above are shown, namely the extreme values, the 3 biggest and the 3 smallest values, surely to be considered as outliers.

First Digits of Final Set   -   {30.0,  17.6,  12.7,  10.1,  8.5,  7.4,  5.7,  4.3,  3.7}
Benford's Law 1st  Digits -  {30.1,  17.6,  12.5,   9.7,  7.9,  6.7,  5.8,  5.1,  4.6}

Figure 5 depicts the digital comparison between this C&F process and Benford's Law. The very low 2.5 SSD value calculated for this result indicates an excellent agreement with the logarithmic, having slightly better fit to Benford as compared with the two previous C&F schemes.

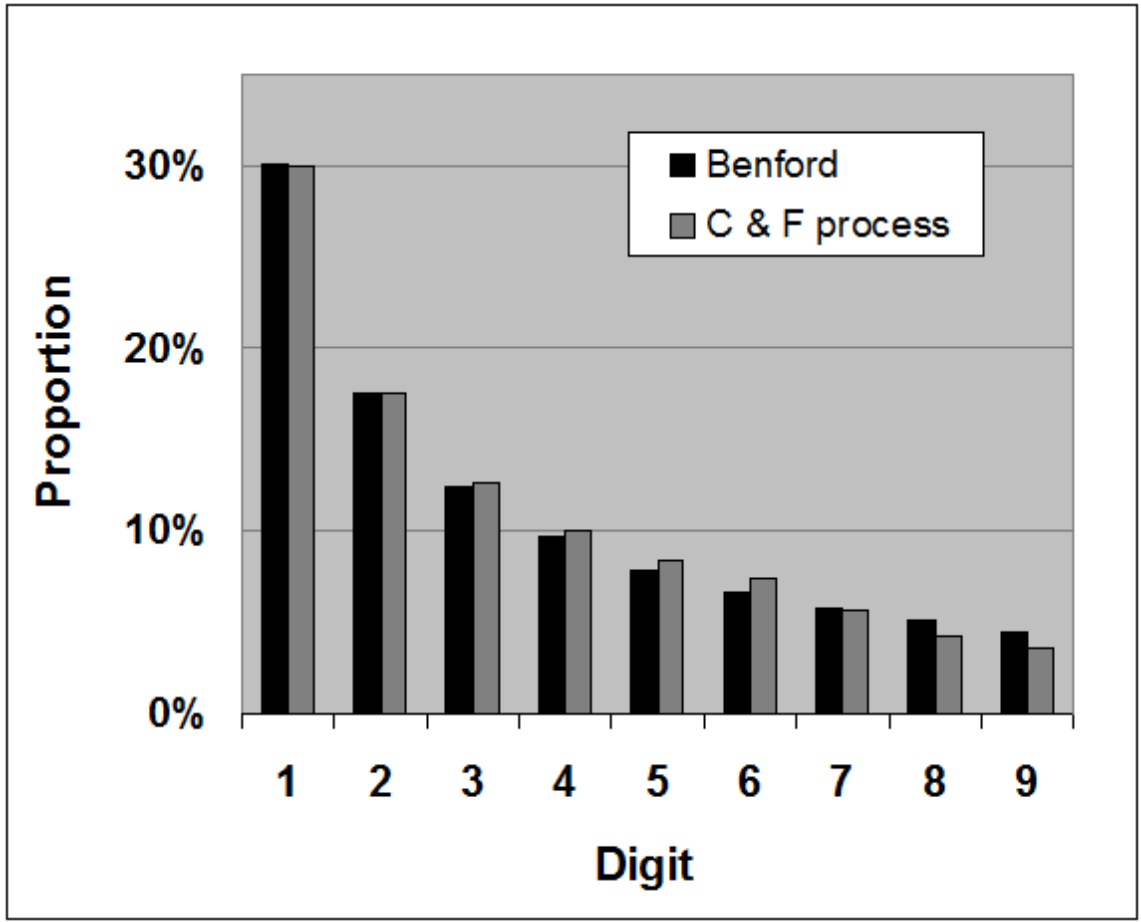

**Figure 5**: C&F Process with a Deterministic Fixed Skewed Ratio 15% - 85%  –   Scheme C

Of all the three C&F schemes above, scheme C yields the closest digital configuration to the Benford proportion. It is superior to scheme B due to its highly skewed split of 15% - 85% which promotes more intense differentiations in resultant relative quantities, and which drives the digital phenomenon. Scheme A also yields at times highly skewed splits near say the 15% - 85% ratio, but very often it's near the more even splits of around 50% - 50%, or other mild splits around say 45% - 55% and 40% - 60%.

Quantitatively, the process takes as input a set of identical numbers with CPOM value of 1, and transforms it into a highly skewed set of numbers where CPOM is $Q_{90\%}/Q_{10\%} = 2.92/0.0018 = 1658$. Here 80% of resultant data falls within (0.0018, 2.92). Figure 6 depicts the histogram from 0 to 8.8 of the final resultant set of numbers after the 3000th binary cycle.

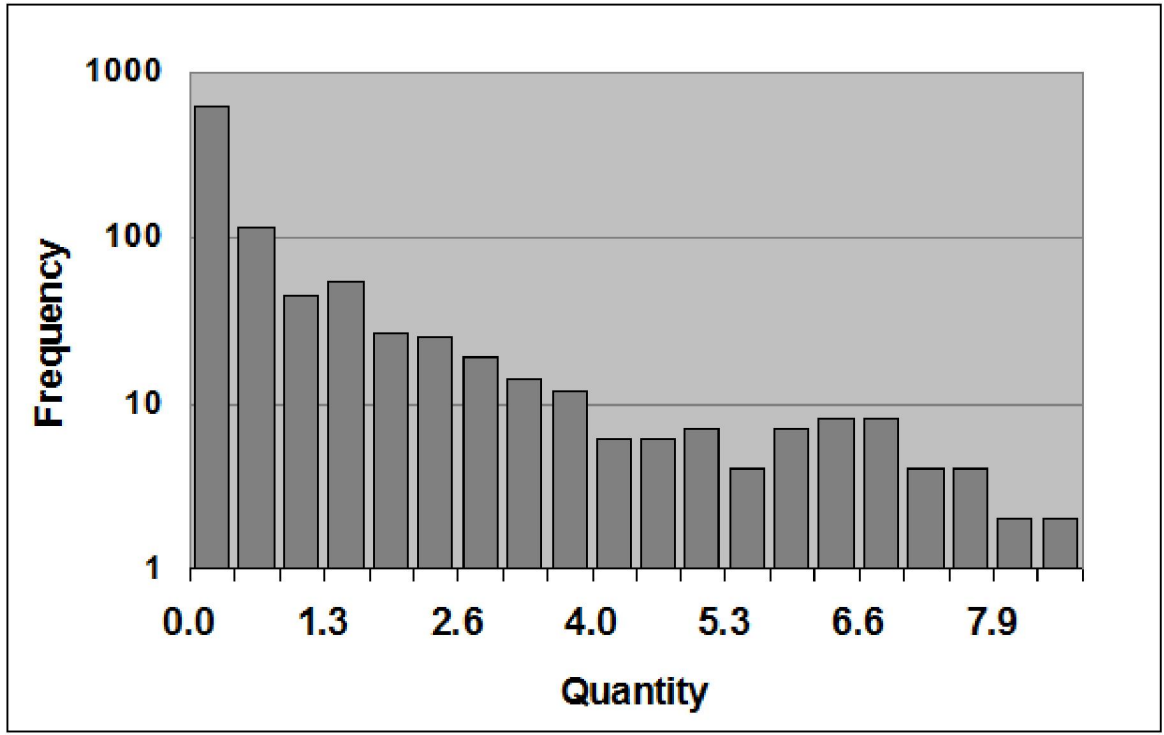

**Figure 6**: Quantitative Configuration with a Fixed Skewed Ratio 15% - 85%  –  Scheme C

As opposed to the case of Random Dependent Partition with sharp differentiation in the rate of convergence to Benford depending on whether the breakup is random or deterministic, here for the Fragmentation and Consolidation model, a deterministic fixed ratio p model yields rapid convergence to Benford just as in the random model. In fact, deterministic models here with highly skewed p and (1 – p) values, such as 15% - 85% for example, seem to converge even faster than random models utilizing the Uniform(0, 1) which is often centered around the more even 50% - 50% or 40% - 60% ratios. The necessary factor driving convergence in C&F models is only the random manner by which the ball which is to be broken is chosen, and the random manner by which the two balls which are to be consolidated are chosen. Since there exists already 'plenty of randomness' in the system in how balls are selected for fragmentations and for consolidations, not much can be gained by adding more randomness with the use of the Uniform(0, 1).

## [3] Examinations of the Random Emergence of the Algebraic Expressions

It should be noted that the first significant digits configuration of the initial set {V, V, V, … L times …V, V, V} is as far removed from the Benford configuration as one could imagine, endowing 0% proportion for 8 digits, and the entire 100% proportion for one privileged digit, namely for that fortunate digit leading V. This fact dramatizes the decisive and very rapid digital transformation taking place under repeated consolidations and fragmentations cycles. Surely, a similar dramatic transformation is occurring at the quantitative level, and actually it is this quantitative transformation which drives the digital one - not the other way around. We start out with a single repeated value V, and a histogram showing a tall and very thin line or a single high bin at V with height of L; then we end up with many diverse and skewed values, represented by a histogram having many bins such as in Figures 2, 4, and 6, and with a tail to the right, where the small is numerous and the big is rare.

Let us examine the mathematical expressions being randomly formed as we record the initial six binary cycles of a C&F process of five balls progressing step-by-step - for one particular scenario of chanced development, with 'authentic' randomness in ball selections intended or imagined by the author. Here the balls split randomly into two parts utilizing the Uniform(0, 1), calling each pair of realization as U and (1 − U), or rather as UJ and (1 − UJ) with J as an index :

{V, V, V, V, V}
{**V**, V, V, V, V} *to be split*
{V*U1, V*(1-U1), V, V, V}
{V*U1, **V*(1-U1)**, V, V, V, **V**} *to be merged*
{V*U1, V*(1-U1) + V, V, V, V}
{**V*U1**, V*(1-U1) + V, V, V, V} to be split
{V*U1*U2, V*U1*(1-U2), V*(1-U1) + V, V, V, V}
{V*U1*U2, V*U1*(1-U2), V*(1-U1) + V, **V**, **V**, V} *to be merged*
{V*U1*U2, V*U1*(1-U2), V*(1-U1) + V, 2V, V}
{V*U1*U2, **V*U1*(1-U2)**, V*(1-U1) + V, 2V, V} *to be split*
{V*U1*U2, V*U1*(1-U2)*U3, V*U1*(1-U2)*(1-U3), V*(1-U1) + V, 2V, V}
{**V*U1*U2**, V*U1*(1-U2)*U3, V*U1*(1-U2)*(1-U3), V*(1-U1) + V, **2V**, V} *to be merged*
{V*U1*U2+2V, V*U1*(1-U2)*U3, V*U1*(1-U2)*(1-U3), V*(1-U1) + V, V}
{V*U1*U2+2V, V*U1*(1-U2)*U3, V*U1*(1-U2)*(1-U3), **V*(1-U1) + V**, V} *to be split*
{V*U1*U2+2V, V*U1*(1-U2)*U3, V*U1*(1-U2)*(1-U3), (V*(1-U1) + V)*U4, (V*(1-U1)+V)*(1-U4), V}
{V*U1*U2+2V, V*U1*(1-U2)*U3, V*U1*(1-U2)*(1-U3), **(V*(1-U1) + V)*U4**, (V*(1-U1)+V)*(1-U4), **V**} *to be merged*
{V*U1*U2+2V, V*U1*(1-U2)*U3, V*U1*(1-U2)*(1-U3), V + (V*(1-U1) + V)*U4, (V*(1-U1) + V)*(1-U4)}
{V*U1*U2+2V, **V*U1*(1-U2)*U3**, V*U1*(1-U2)*(1-U3), V + (V*(1-U1) + V)*U4, (V*(1-U1) + V)*(1-U4)} *to be split*
{V*U1*U2+2V, V*U1*(1-U2)*U3*U5, V*U1*(1-U2)*U3*(1-U5), V*U1*(1-U2)*(1-U3), V + (V*(1-U1) + V)*U4, (V*(1-U1) + V)*(1-U4)}
{V*U1*U2+2V, V*U1*(1-U2)*U3*U5, **V*U1*(1-U2)*U3*(1-U5)**, V*U1*(1-U2)*(1-U3), **V+(V*(1-U1) + V)*U4**, (V*(1-U1) + V)*(1-U4)} *to be merged*
{V*U1*U2+2V, V*U1*(1-U2)*U3*U5, V*U1*(1-U2)*U3*(1-U5) + V+(V*(1-U1) + V)*U4, V*U1*(1-U2)*(1-U3), (V*(1-U1) + V)*(1-U4)}
{**V*U1*U2+2V**, V*U1*(1-U2)*U3*U5, V*U1*(1-U2)*U3*(1-U5) + V+(V*(1-U1) + V)*U4, V*U1*(1-U2)*(1-U3), (V*(1-U1) + V)*(1-U4)} *to be split*
{(V*U1*U2+2V)*U6, (V*U1*U2+2V)*(1-U6), V*U1*(1-U2)*U3*U5, V*U1*(1-U2)*U3*(1-U5) + V+(V*(1-U1) + V)*U4, V*U1*(1-U2)*(1-U3), (V*(1-U1) + V)*(1-U4)}
{(V*U1*U2+2V)*U6, (V*U1*U2+2V)*(1-U6), **V*U1*(1-U2)*U3*U5**, V*U1*(1-U2)*U3*(1-U5) + V+(V*(1-U1) + V)*U4, V*U1*(1-U2)*(1-U3), **(V*(1-U1) + V)*(1-U4)**} *to merge*
{(V*U1*U2+2V)*U6, (V*U1*U2+2V)*(1-U6), V*U1*(1-U2)*U3*U5 + (V*(1-U1) + V)*(1-U4), V*U1*(1-U2)*U3*(1-U5) + V+(V*(1-U1) + V)*U4, V*U1*(1-U2)*(1-U3)}

**24**

Let's examine the five terms after the last cycle:

(V*U1*U2 + 2V)*U6
(V*U1*U2 + 2V)*(1 - U6)
V*U1*(1 - U2)*U3*U5 + (V*(1 - U1) + V)*(1 - U4)
V*U1*(1 - U2)*U3*(1 - U5) + V + (V*(1 - U1) + V)*U4
V*U1*(1 - U2)*(1 - U3)

The initial V weight of each ball appears as a singular factor only once, since in essence it represents the scale of the entire system. Dividing by V (or equivalently setting V = 1) we obtain the 'pure' or dimensionless set of algebraic expressions:

(U1*U2 + 2)*U6
(U1*U2 + 2)*(1 - U6)
U1*(1 - U2)*U3*U5 + ((1 - U1) + 1)*(1 - U4)
U1*(1 - U2)*U3*(1 - U5) + 1+ ((1 - U1) + 1)*U4
U1*(1 - U2)*(1 - U3)

Simplifying a bit, we get:

(U1*U2 + 2)*U6
(U1*U2 + 2)*(1 - U6)
U1*U3*U5*(1 - U2) + (2 - U1)*(1 - U4)
U1*U3*(1 - U2)*(1 - U5) + 1+ (2 - U1)*U4
U1*(1 - U2)*(1 - U3)

**It is necessary to keep in mind the randomness involved in building up and writing these five algebraic expressions above.** In any case, for this particular random trajectory of events above, out of five expressions, two are sums (minority), and three are products (majority). In fact, since each full cycle yields **3** newly created balls, **2** of which are multiplicative, namely X(p) and X*(1 – p), and **1** of which is an additive, therefore the statistical tendency of the system after numerous cycles is to have approximately 2/3 multiplicative expressions and 1/3 additive expressions. The reason this is only approximately so is that a term such as X(p) converts X into a product if X itself was additive, but it leaves it as a product if X itself was already a product. In the same vein, a term such as X + Y converts X and Y into a sum if X and Y themselves were products, but it leaves them as a sum if X and Y themselves were already sums.

On the face of it, the existence of additive expressions in about one third of all expressions, does not bode very well for Benford digital configuration; the glass is two-third full and one-third empty, and one doesn't know whether he or she should be happy or should be sad. Clearly, the manifestation of tugs of war that occur between addition and multiplication in all Fragmentation and Consolidation models needs to be investigated further.

An abstract Monte Carlo computer simulation program of C&F scheme utilizing a deterministic fixed p and (1 – p) ratios was run, with 43 identical initial values of V, experiencing 52 full cycles. The program was designed to obtain the random algebraic expressions of the process, but not to perform any numerical calculation, leaving V and p as variables. Additive terms are shown as bold types for emphasis. The results after the 52nd cycle are as follows:

1) ((((V)p)p+(((V)p)(1-p))p)p)p
2) ((V)(1-p))p
3) ((V)p)p **+** (V)
4) (V+((V)p)(1-p)+(V)p+(V)p+(V)(1-p))p
5) (V+((V)p)p)(1-p)p
6) V **+** V **+** ((V)p)p **+** V
7) (V+(V)(1-p)+((V)p)(1-p) **+** ((V+V)p+V)(1-p) **+** ((((V)p)p+(((V)p)(1-p))p)(1-p))p)p
8) (V)p **+** (V)(p)
9) (V)(1-p)
10) (V+(V)(1-p)+V)(1-p)
11) (V+((V)p)p)p **+** ((V)(1-p)+V)p
12) (V)p
13) ((V)(1-p))p **+** V **+** V **+** (((V)p)p)(1-p) **+** V
14) (((V)p)p)p **+** (V)(1-p) **+** (V)p
15) ((V+(V)(1-p)+V)p)(1-p)
16) ((((V)p)(1-p))(1-p))(1-p) **+** ((((V)p)p+(((V)p)(1-p))p)(1-p))(1-p) **+** (V)p
17) (((V+V)p)(1-p))p
18) (V+V)(1-p)
19) (V+((V)p)p)(1-p)(1-p)
20) (V)(1-p) **+** V **+** (V+((V)p)(1-p)+(V)p+(V)p+(V)(1-p))(1-p)
21) (V+(V)(1-p)+((V)p)(1-p)+((V+V)p+V)(1-p)+((((V)p)p+(((V)p)(1-p))p)(1-p))p)(1-p)
22) (V)p **+** (V)(1-p)
23) (((V)p)(1-p))p
24) ((((V+V)p)(1-p))(1-p))p
25) ((V)p)(1-p) **+** (V)p
26) (V)p **+** ((V)p)p **+** ((V)(1-p))(1-p)
27) ((V+(V)(1-p)+V)p)p
28) ((V)(1-p))p
29) ((((V)p)p+(((V)p)(1-p))p)p)(1-p)
30) ((V+V)p+V)p **+** V **+** (V)(1-p) **+** (V)(1-p) **+** ((V+V)p)p
31) (((((V)p)(1-p))(1-p))p)p
32) (V)(1-p) **+** (V)(1-p)
33) ((V)p)(1-p)
34) ((V)p+V)p
35) (V+V)(1-p) **+** (V)(1-p)
36) (((V)p)(1-p))(1-p)
37) ((V)(1-p)+V)(1-p) **+** (V)(1-p) **+** (V)(1-p) **+** ((V)p)(1-p) **+** V **+** V **+** ((V)p)p
38) ((((V+V)p)(1-p))(1-p))(1-p)
39) (V)p **+** ((V)(1-p))(1-p) **+** ((V)p+V)(1-p)
40) ((V)(1-p))(1-p)
41) ((V)(1-p))p
42) ((V)(1-p))(1-p) **+** V
43) (((((V)p)(1-p))(1-p))p)(1-p)

Out of 43 final expressions at the end of the 52th cycle, 18 are additive and 25 are multiplicative. Simulation runs with 1000+ balls and with 4000+ full cycles consistency show nearly exact ratios of 2/3 multiplicative expressions and 1/3 additive expressions at the end of the processes.

## [4] Tugs of War between Addition and Multiplication in C&F Models

Tugs of war between additions and multiplications in the context of Benford's Law are discussed in an article by this author titled "Arithmetical Tugs of War and Benford's Law" on https://arxiv.org/. Link to that article can be found at: https://arxiv.org/abs/1410.2174

Let us summarize the main results from that article that are relevant to C&F processes:

Random multiplication processes favor the small over the big - leading to skewed data.
Random multiplication processes induce two essential results:
(A) A dramatic increase in skewness – an essential criterion for Benford behavior.
(B) An increase in the order of magnitude – another essential criterion for Benford behavior.
Hence multiplication processes are highly conducive to Benford behavior.

Random addition processes favor the medium over the small and over the big.
Random addition processes do not induce any results that are essential in the criterion for Benford behavior:
(A) Lacking any increase in skewness, and even actively increasing the symmetry of resultant distribution, with added concentration forming around the center/medium.
(B) Lacking any increase in order of magnitude beyond the existing maximum order of magnitude within the set of added variables.
Hence addition processes are highly detrimental to Benford behavior.

<u>Addition Processes:</u>     Less POM  - More Symmetry -  CLT  -  Normal  - Anti Benford
<u>Multiplication Processes:</u>  More POM  - More Skewness -  MCLT - Lognormal - Pro Benford

In addition, it is necessary to point out to a significant limitation in the effectiveness of the Central Limit Theorem. The CLT's Achilles' heel - in terms of its rate of convergence to the Normal - is the possibility that added variables are highly skewed and that they come with very high order of magnitude. This is a bad combination for the CLT. Except for Uniforms, Normals, and other symmetrical distributions which converge to the Normal quite quickly after very few additions regardless of the value of OOM of added variables, all other asymmetrical (skewed) distributions show a distinct rate of convergence depending on their OOM value. For skewed variables, whenever OOM is of very high value, CLT can manifest itself with difficulties, and very slowly, only after a truly large number of additions of these random variables. On the other hand, when skewed variables are of very low OOM value, CLT achieves near Normality quite quickly after only very few additions.

For the Consolidations and Fragmentation model, each fragmentation process contributes to the system two products - each with a minimum of 2 multiplicands and possibly more; similarly each consolidation process contributes to the system an additive expression with a minimum of 2 addends and possibly more; and therefore there exists a tug of war here between additions and multiplications with respect to Benford behavior; a struggle between the Central Limit Theorem and the Multiplicative Central Limit Theorem; between the Normal and the Lognormal. Remarkably, even though Benford frequently loses numerous C&F battles, yet he wins the war in the long run. One existing feature here that is partially saving the system from deviation from Benford is that on average only about one-third of the expressions are additive; and that even within those expressions there are plenty of arithmetical multiplicative elements involved. Surely there are some additive expressions with 3 addends which might be quite detrimental to Benford, but they are far and few between; and there are even more menacing expressions with 4 addends, but luckily these are even rarer.

The general [theoretical] understanding gained in that article regarding multiplication and addition processes enables us to thoroughly explain the [empirical] strong logarithmic behavior in all Consolidations and Fragmentations models; namely the reason addition effects do not manage to significantly retard multiplication effects. In a nutshell, the C&F process is Benford because it uses the high OOM variable of Uniform(0, 1) which contributes to high order of magnitude, skewness, and thus Benfordness. As a consequence, the C&F process encounters the Achilles' heel of the Central Limit Theorem and additions are not very effective. Order of magnitude of the Uniform(0, 1) calculated as LOG(1/0) is infinite. CPOM calculated as $Q_{90\%}$ divided by $Q_{10\%}$ is 0.9/0.1 or simply 9. Surely, the C&F model cannot use any low OOM variable such as say Uniform(5, 7), because it needs to break a whole quantity into two fractions, and this can only be achieved via the high OOM variable Uniform(0, 1). Such high OOM values, coupled with the fact that the terms within the additive expressions almost always involve also some multiplications (which are always skewed), guarantee that the Central Limit Theorem is very slow to act here and that its retarded rate of convergence does not manage to even begin to ruin the general multiplicative tendencies of the system. Skewness for these multiplicative terms hiding within the additive expressions is guaranteed by virtue of simply being multiplication. All multiplication processes yield skewed set of values. Since the vast majority of the additive terms here are with only 2, 3, or 4 addends, the Central Limit Theorem does not even begin to manifests itself.

And what about deterministic models of fixed 50% - 50% ratios or fixed 15% - 85% ratios, etc.?! Even though these deterministic models do not apply the Uniform(0, 1) in any way, the effective results of such deterministic p and (1 – p) ratios are also of high OOM values. Occasional terms such as V*p*p*(1-p)*(1-p), or V*p*p*p*p, and so forth, are of very low value since p < 1, while a fortunate surviving V without any p factor stays as large as V itself, and all that implies high variability. Even a comparison of V*p with V*p*p*p*p indicates the possibility of sufficiency high OOM in the system since the value of the former is much higher than the value of the latter. For example, even for the fixed balanced 50% - 50% ratios model, variations between possible terms are very high, such as say between V*0.5 and V*0.5*0.5*0.5*0.5 = V*0.0625, guaranteeing sufficiently high order of magnitude and preventing the Central Limit Theorem from manifesting itself.

For both models, for the random utilizing Uniform(0, 1), and for the deterministic utilizing fixed p and (1 – p) ratios, order of magnitude may not be substantial, but it is sufficient to significantly retard CLT given that CLT has only 2 or 3 addends to work with.

## [5] Examination of C&F Processes Via Gradual Run Reveals Saturation

Another run of scheme A is performed in order to demonstrate its rate of convergence to the logarithmic, as well as to demonstrate the digital saturation that occurs after sufficient number of cycles. The scheme starts out with 2000 balls having an initial identical value of 1, and utilizing the Uniform(0, 1) for fragmentations, but instead of running 8000 binary cycles and stopping, here it is extended to 13000 cycles, and the scheme is executed slowly and gradually, in stages, in order to be able to take snapshots occasionally, and to examine how first digits, SSD, and CPOM gradually evolve. Figure 7 depicts the computerized results in details. Evidently, it takes the system about 2*L cycles [*namely 2*2000 = 4000 cycles*] to achieve its nearly perfect Benford behavior. Beyond 4000 cycles, nothing new is achieved, although the system vigilantly maintains and guards the logarithmic status that it has thus obtained. CPOM though seems to be growing continuously well beyond these 4000 cycles, although it is possible that it might change course and reverse itself further on perhaps.

| # of Cycles | 1 | 2 | 3 | 4 | 5 | 6 | 7 | 8 | 9 | SSD | $Q_{10}$ | $Q_{90}$ | CPOM |
|---|---|---|---|---|---|---|---|---|---|---|---|---|---|
| 0 | 100% | 0% | 0% | 0% | 0% | 0% | 0% | 0% | 0% | 5634 | 1.00000 | 1.00 | 1 |
| 500 | 59% | 14% | 6% | 5% | 4% | 4% | 3% | 3% | 3% | 951 | 0.20131 | 2.00 | 10 |
| 1000 | 45% | 17% | 8% | 7% | 5% | 5% | 5% | 4% | 3% | 253 | 0.09119 | 2.00 | 22 |
| 1500 | 38% | 18% | 11% | 8% | 6% | 5% | 5% | 5% | 3% | 79 | 0.04417 | 2.29 | 52 |
| 2000 | 35% | 19% | 12% | 7% | 7% | 6% | 5% | 4% | 4% | 38 | 0.02965 | 2.47 | 83 |
| 2500 | 33% | 19% | 12% | 8% | 7% | 6% | 5% | 5% | 4% | 18 | 0.02116 | 2.64 | 125 |
| 3000 | 33% | 19% | 12% | 9% | 7% | 6% | 5% | 5% | 4% | 14 | 0.01410 | 2.69 | 191 |
| 4000 | 31% | 17% | 13% | 10% | 8% | 6% | 5% | 5% | 5% | 2 | 0.00915 | 2.84 | 311 |
| 5000 | 32% | 18% | 13% | 11% | 7% | 7% | 5% | 5% | 5% | 6 | 0.00346 | 2.96 | 856 |
| 6000 | 31% | 18% | 12% | 11% | 7% | 6% | 5% | 5% | 4% | 5 | 0.00296 | 2.97 | 1003 |
| 7000 | 31% | 18% | 12% | 10% | 8% | 7% | 5% | 5% | 4% | 2 | 0.00188 | 2.97 | 1581 |
| 8000 | 30% | 18% | 13% | 11% | 7% | 6% | 6% | 5% | 4% | 2 | 0.00121 | 3.09 | 2543 |
| 10000 | 29% | 19% | 13% | 10% | 8% | 5% | 6% | 5% | 5% | 5 | 0.00093 | 2.99 | 3209 |
| 13000 | 31% | 16% | 12% | 10% | 8% | 6% | 5% | 7% | 4% | 6 | 0.00086 | 2.94 | 3408 |

**Figure 7**: The Evolution of First Digits, SSD, and CPOM of a Gradual Run of Scheme A

## [6] Related Log Conjecture Justifies Logarithmic Behavior of C&F Models

Examining the histogram of log values of gradual run Scheme A - helps in relating the result of C&F models to Related Log Conjecture discussed in Kossovsky (2014) chapter 63. Figure 8 depicts this histogram after 8000 cycles. Since the original V value for all the balls was 1, it stands to reason that even after 8000 cycles the histogram is at its highest for values around this point which corresponds to the log value of LOG(1) = 0. Earlier in Figure 7 it was shown that after 8000 cycles, $Q_{10}$ = 0.00121, and $Q_{90}$ = 3.09, corresponding to the log values of LOG(0.00121) = -2.92 and LOG(3.09) = +0.49. Indeed, one can easily ascertain that in Figure 8 the core of the data (the central 80% of the data) is approximately between -3 and +0.5. All the requirements of Related Log Conjecture are easily satisfied for the curve in Figure 8, especially the requirement of having wide enough span on the log-axis, being at least 5 units here, and which is comfortably more than the usual 3 or 4 units required for strong logarithmic behavior. In conclusion: Benford behavior here is in perfect harmony and consistent with Related Log Conjecture.

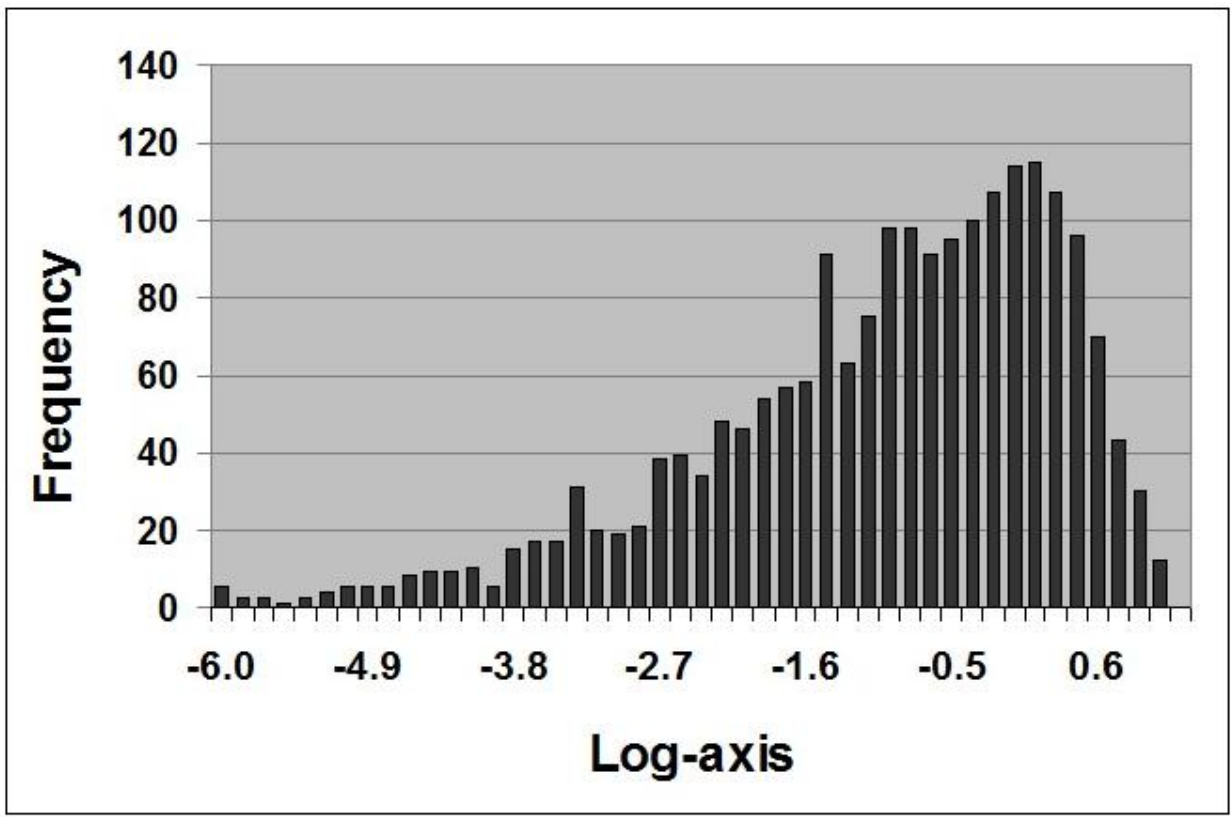

**Figure 8**: Log Histogram Complies with Related Log Conjecture – Scheme A, 8000 Cycles

# [7] Partial Convergence to Benford for Models with Few Initial Quantities

Another C&F scheme is run via Monte Carlo computer simulations. The scheme starts out with only 100 balls having an initial identical value of 1, and utilizing the Uniform(0, 1) for fragmentations. It is executed slowly and gradually in stages for a total of 9000 cycles, in order to be able to take snapshots occasionally, and to examine how first digits, SSD, and CPOM gradually evolve. Surely such a meager set of only 100 balls does not have sufficient number of quantities to fully converge to Benford, although the system achieves partial convergence in spite of its small size.

Figure 9 depicts the computerized results in details. Evidently, after about 2*L to 3*L cycles (2*100 to 3*100), namely 200 to 300 cycles, the system achieves some kind of digital stability with its convergence to an approximate Benford configuration with SSD stabilizing around the value of about 50. Beyond about 200 or 300 cycles, nothing new is achieved, although the system maintains its approximate logarithmic status it has thus obtained. CPOM on the other hand, grows steadily initially for the first 300 cycles or so; but then it starts reversing itself, fluctuating higher and lower as more cycles are added.

| # of Cycles | 1 | 2 | 3 | 4 | 5 | 6 | 7 | 8 | 9 | SSD | $Q_{10}$ | $Q_{90}$ | CPOM |
|---|---|---|---|---|---|---|---|---|---|---|---|---|---|
| 0 | 100% | 0% | 0% | 0% | 0% | 0% | 0% | 0% | 0% | 5634 | 1.00000 | 1.00 | 1 |
| 50 | 42% | 16% | 9% | 7% | 4% | 4% | 6% | 7% | 5% | 190 | 0.07469 | 2.05 | 27 |
| 100 | 31% | 19% | 15% | 5% | 5% | 7% | 4% | 3% | 11% | 89 | 0.04178 | 2.53 | 61 |
| 150 | 28% | 12% | 13% | 13% | 7% | 7% | 9% | 5% | 6% | 60 | 0.03510 | 2.55 | 73 |
| 200 | 28% | 12% | 18% | 10% | 5% | 9% | 10% | 6% | 2% | 105 | 0.01072 | 3.11 | 290 |
| 250 | 21% | 16% | 15% | 12% | 9% | 8% | 9% | 8% | 2% | 125 | 0.00567 | 3.23 | 570 |
| 300 | 25% | 21% | 13% | 9% | 6% | 8% | 5% | 8% | 5% | 53 | 0.00156 | 2.79 | 1790 |
| 350 | 27% | 22% | 9% | 14% | 4% | 8% | 4% | 7% | 5% | 84 | 0.00262 | 2.77 | 1059 |
| 800 | 36% | 20% | 8% | 9% | 5% | 7% | 4% | 6% | 5% | 74 | 0.00034 | 2.14 | 6296 |
| 2000 | 28% | 16% | 17% | 7% | 9% | 7% | 8% | 7% | 1% | 57 | 0.00231 | 3.18 | 1377 |
| 5000 | 26% | 20% | 10% | 12% | 7% | 8% | 7% | 6% | 4% | 39 | 0.00079 | 3.03 | 3853 |
| 9000 | 29% | 21% | 11% | 7% | 12% | 8% | 4% | 3% | 5% | 48 | 0.00087 | 1.91 | 2182 |

**Figure 9**: First Digits, SSD, and CPOM of a Gradual Run Scheme with Only 100 Balls

Jan 2, 2019
Alex Ely Kossovsky
akossovsky@gmail.com